\documentclass[15pt, twoside]{article}
\usepackage{amssymb}
\usepackage[all]{xy}
\usepackage{amsthm,amsmath,color}
\usepackage{pgfmath}
\usepackage{tikz-cd}
\usetikzlibrary{arrows, matrix}
\usepackage{zref-abspage}
\usepackage{perpage}
\usepackage[top=2.5cm,right=3.5cm,bottom=2.5cm,left=2.5cm]{geometry}
\usepackage{graphicx}
\usepackage{fancyhdr}
\usepackage{makeidx}
\usepackage{ifpdf}
\usepackage[colorlinks=true, linkcolor=blue, linktoc=all, citecolor=magenta, filecolor=cyan, urlcolor=cyan, linkbordercolor={1 0 0}, citebordercolor={0 1 0}, urlbordercolor={0 1 1}]{hyperref}
\newtheorem{theorem}{Theorem}[section]

\newtheorem{proposition}{Proposition}[section]

\newtheorem{example}{Example}[section]
\begin{document}
\fancyhead{}
\fancyfoot{}
\chead[M. J. Afshari, S. Varsaie]{Universal Super Vector Bundles}
\lhead[\thepage]{}
\rhead[]{\thepage}
\pagestyle{fancy}
\author{M. J. Afshari, S. Varsaie\footnote{\textit{E-mail addresses:} \href{mailto: afshari.mj@iasbs.ac.ir}%
{afshari.mj@iasbs.ac.ir} (M. J. Afshari), \href{mailto: varsaie@iasbs.ac.ir}%
{varsaie@iasbs.ac.ir} (S. Varsaie).}
\\
\small
\textit{Department of Mathematics,}
\\
\small
\textit{Institute for Advanced Studies in Basic Sciences (IASBS),}
\\
\small
\textit{Zanjan 45137-66731, Iran}}                            
\title{Universal Super Vector Bundles}
\renewcommand{\baselinestretch}{1.5}
\setlength{\baselineskip}{1.5\baselineskip}
\setcounter{page}{1}
\maketitle
\begin{abstract}
A new generalization of Grassmannians, called $\nu$-grassmannians, and a canonical super vector bundle over this new space, say $\Gamma$, are introduced. Then, constructing a Gauss supermap of a super vector bundle, the universal property of $\Gamma$ is discussed. Finally, we generalize one of the main theorems of homotopy classification for vector bundles in supergeometry.
\end{abstract}
\section{Introduction}
This paper aims at deriving a homotopy classification for super vector bundles. The importance of it is in connection with finding proper generalization of Chern classes in supergeometry. Indeed, Chern classes are some cohomology elements associated to isomorphism classes of complex vector bundles in common geometry. In the category of vector bundles,  it is shown that canonical vector bundles $\gamma^n_k$ on grassmannians $Gr(n,k)$ are universal. Equivalently, associated to each vector bundle $\cal{E}$ on M, up to homotopy, there exists a unique map $f:M\to Gr(n,k)$, for sufficiently large $n$, such that $\cal{E}$ is isomorphic with the induced bundle of $\gamma^n_k$ under $f$. In addition, Chern classes of a vector bundle may be described as the pullback of Chern classes of the universal bundle.
\\
To have an appropriate generalization of homotopy classification theorem, one should have a proper generalization of grassmannian. Supergrassmannian, introduced in \cite{Manin} and Grassmannian, in some sense, are homotopy equivalent (cf. subsection 2.2). Therefore, cohomology group associated to supergrassmannian is equal to that of grassmannian. In other words, the former group contains no information about superstructure. Hence, from classifying space viewpoint, supergrassmannians are not good generalization of Grassmannians.
\\
In this paper, first, following \cite{Bahadory}, we introduce $\nu$-grassmannians denoted by $_{\nu}Gr(m|n)$, as a new supergeneralization of Grassmannians. In addition, we show the existence of $\Gamma$, a canonical super vector bundle over $_{\nu}Gr(m|n)$. After introducing Gauss supermaps for super vector bundles, universal property of $\Gamma$ is studied. At the end, we extend one of the main theorems on homotopy classification for vetor bundles to supergeometry.
\\
 There are different approaches to generalize  Chern classes in supergeometry, such as homotopy or analytic approach. In this paper our approach is homotopic. Although there are not many articles with homotopy approach, but one may refer to \cite{V-manin} as a good example for such papers. Nevertheless, much more efforts have been made for generalizing Chern classes in supergeometry by analytic approach. One may refer to \cite{Bartocci-B}, \cite{Bartocci-B-H}, \cite{B-H}, \cite{Landi}, \cite{Manin-P}, \cite{V-Manin-P}. But, in all these works, the classes obtained in this way are nothing but the Chern classes of the reduced vector bundle(s) and they do not have any information about the superstructure.  
\section{Preliminaries}
In this section, first, we recall some basic definitions of supergeometry. Then, we introduce a supergeneralization of Grassmannian called $\nu$-grassmannian.
\subsection{Supermanifolds}
A \textit{super ringed space} is a pair $(X, \mathcal{O})$ where $X$ is a topological space and $\mathcal{O}$ is a sheaf of commutative $\mathbb{Z}_{2}$-graded rings with units on $X$. Let $\mathcal{O}(U)= \mathcal{O}^{ev}(U) \oplus \mathcal{O}^{odd}(U)$ for any open subset $U$ of $X$. An element $a$ of $\mathcal{O}(U)$ is called a homogeneous element of parity $p(a)=0$ if $a \in \mathcal{O}^{ev}(U)$  and it is a homogeneous element of parity $p(a)=1$ if $a \in \mathcal{O}^{odd}(U)$. A morphism between two super ringed spaces $(X, \mathcal{O}_{X})$ and $(Y, \mathcal{O}_{Y})$ is a pair $(\widetilde{\psi}, \psi^{*})$ such that $\widetilde{\psi}: X \longrightarrow Y$ is a continuous map and $\psi^{*}: \mathcal{O}_{Y} \longrightarrow \widetilde{\psi}_{*}(\mathcal{O}_{X})$ is a homomorphism between the sheaves of $\mathbb{Z}_{2}$-graded rings.
\\
A \textit{superdomain} is a super ringed space $\mathbb{R}^{p|q}:=(\mathbb{R}^{p}, \mathcal{O})$ where
\begin{equation*}
\mathcal{O}= \mathbf{C}^{\infty}_{\mathbb{R}^{p}} \otimes_{\mathbb{R}} \wedge \mathbb{R}^{q}, \qquad p, q \in \mathbb{N}.
\end{equation*}
By $\mathbf{C}^{\infty}_{\mathbb{R}^{p}}$ we mean the sheaf of smooth functions on $\mathbb{R}^{p}$. A super ringed space which is locally isomorphic to $\mathbb{R}^{p|q}$ is called a \textit{supermanifold} of dimension $p|q$. Note that a morphism $(\widetilde{\psi}, \psi^{*})$ between two supermanifolds $(X, \mathcal{O}_{X})$ and $(Y, \mathcal{O}_{Y})$ is just a morphism between the super ringed spaces such that for any $x \in X$, $\psi^*: \mathcal{O}_{Y, \widetilde{\psi}(x)} \longrightarrow \widetilde{\psi}_{*}(\mathcal{O}_{X, x})$ is local, i.e., $\psi^{*}(m_{\widetilde{\psi}(x)}) \subseteq m_{x}$, where $m_x$ is the unique maximal ideal in $\mathcal{O}_{X, x}$.

\subsection{$\nu$-grassmannian}\label{grassmannian}
Supergrassmannians are not good generalization of Grassmannians. Indeed these two, in some sense, are homotopy equivalent. This equivalency may be shown easily in the case of projective superspaces.
\\
To this end, let $\mathbb{P}^{m|n} = (\mathbb{RP}^{m}, \mathcal{O}_{m|n})$ be the real projective superspace.
By a \textit{deformation retraction} from $\mathcal{O}_{m|n}$ to $\mathcal{O}_{m}$, the sheaf of rings of all real valued functions on the real projective space $\mathbb{P}^{m}$, we mean that there is a sheaf of $\mathbb{Z}_{2}$-graded rings, say $\mathcal{A}$, on $\mathbb{RP}^{m} \times \mathbb{R}$ such that $\mathcal{A}_{\mathbb{RP}^{m} \times \{0\}} = \mathcal{O}_{m|n}$, and $\mathcal{A}_{\mathbb{RP}^{m} \times \{1\}} = \mathcal{O}_{m|n}$. In addition, there are morphisms $r: \mathcal{O}_{m} \longrightarrow \mathcal{O}_{m|n}$, $j: \mathcal{O}_{m|n} \longrightarrow \mathcal{O}_{m}$ and $H: \mathcal{O}_{m|n} \longrightarrow \mathcal{A}$ along $\mathbb{RP}^{m} \times \mathbb{R} \longrightarrow \mathbb{RP}^{m}$ with $(z, t) \longmapsto z$ such that $j_{0} \circ H=id$ and $j_{1} \circ H= r \circ j$ where $j_{0}: \mathcal{A} \longrightarrow \mathcal{O}_{m|n}$ and $j_{1}: \mathcal{A} \longrightarrow \mathcal{O}_{m|n}$ are morphisms of sheaves respectively satisfying the following conditions:
\begin{equation*}
x_{i} \otimes 1 \longmapsto x_{i}, \qquad e_{j} \otimes 1 \longmapsto e_{j}, \qquad 1 \otimes t \longmapsto 0,
\end{equation*}
\begin{equation*}
x_{i} \otimes 1 \longmapsto x_{i}, \qquad e_{j} \otimes 1 \longmapsto e_{j}, \qquad 1 \otimes t \longmapsto 1.
\end{equation*}
\begin{proposition}
There is a deformation retraction from $\mathcal{O}_{m|n}$ to $\mathcal{O}_{m}$.
\end{proposition}
\textit{Proof}:
Set $\mathcal{A} = \mathcal{O}_{m|n} \otimes_{\mathbb{R}} \mathbf{C}^{\infty}_{\mathbb{R}}$ and $U_i=\mathbb{R}^m$, $1 \leq i \leq m+1$. Now, consider $H: \mathcal{O}_{m|n} \longrightarrow \mathcal{A}$ which is locally defined along $U_{i} \times \mathbb{R} \longrightarrow U_{i}$ as follows:
\begin{equation*}
x_{k} \longmapsto x_{k}, \qquad e_{l} \longmapsto (1-t)e_{l}.
\end{equation*}
One may easily show that $j_{0} \circ H=id$ and $j_{1} \circ H= r \circ j$ where $r: \mathcal{O}_{m} \longrightarrow \mathcal{O}_{m|n}$ is a morphism with $x_{k} \longmapsto x_{k}$ and $j: \mathcal{O}_{m|n} \longrightarrow \mathcal{O}_{m}$ is a morphism with $x_{k} \longmapsto x_{k}$ and $e_{l} \longmapsto 0$.
\begin{flushright}
$\square$
\end{flushright}
This proposition shows that in the construction of projective superspaces, the odd variables do not play principal roles. Solving this problem
is our motivation for defining $\nu$- Projective spaces or generally $\nu$- grassmannians. Before that, it is necessary to recall some basic concepts.
\\
A $\nu$-\textit{domain} of dimension $p|q$ is a super ringed space $\mathbb{R}^{p|q}$ which carries an odd involution $\nu$, i.e.,
\begin{equation*}
\nu: \mathcal{O} \longrightarrow \mathcal{O}, \qquad \nu(\mathcal{O}^{ev}) \subseteq \mathcal{O}^{odd}, \qquad \nu(\mathcal{O}^{odd}) \subseteq \mathcal{O}^{ev}, \qquad \nu^{2}=id.
\end{equation*}
In addition, $\nu$ is a homomorphism between $\mathbb{C}^{\infty}$-modules.
\\
Let $k$, $l$, $m$ and $n$ be non-negative integers with $k<m$ and $l<n$. For convenience from now on, we set $p=k(m-k)+l(n-l)$ and $q=k(n-l)+l(m-k)$.
\\
A \textit{real} $\nu$-\textit{grassmannian}, $_{\nu}Gr_{k|l}(m|n)$, or shortly $_{\nu}Gr= (Gr_{k}^{m} \times Gr_{l}^{n}, \mathcal{G})$, is a real superspace obtained by gluing $\nu$-domains $(\mathbb{R}^{p}, \mathcal{O})$ of dimension $p|q$.
\\
Here, we need to set some notations that are useful later.
\\
Let $I$ be a $k|l$ multi-index, i.e., an increasing finite sequence of $\{1, \cdots, m+n\}$ with $k+l$ elements. So one may set
\begin{equation*}
I:= \{i_{a}\}_{a=1}^{k+l}.
\end{equation*}
A standard $(k|l) \times (m|n)$ supermatrix, say $T$, may be decomposed into four blocks as follows:
$$\left[ \! \! \!
\begin{tabular}{ccc}
$A_{k\times m}$   &  $\vline$  & $B_{k \times n}$\\
---  ---  ---             &  $\vline$  & ---  ---  ---          \\
$C_{l \times m}$   &   $\vline$  & $D_{l \times n}$
\end{tabular}
\! \! \! \right]$$
The upper left and lower right blocks are filled by even elements. These two blocks together form the even part of $T$. The upper right and lower left blocks are filled by odd elements and form the odd part of $T$.
\\
The columns with indices in $I$ together form a minor denoted by $M_I(T)$. 
\\
A pseudo-unit matrix $id_{I}$ corresponding to $k|l$ multi-index $I$, is a $(k|l) \times (k|l)$ matrix whose all entries are zero except on its main diagonal that are $1$ or $1\nu$, where $1\nu$ is a formal expression used as odd unit. For each open subset $U$ of $\mathbb{R}^p$ and each $z \in \mathcal{O}(U)$, we also need the following rules:
\begin{equation*}
z.(1\nu):= \nu(z), \qquad (1\nu)(1\nu)=1.
\end{equation*}
So for each $I$, one has
\begin{equation*}
id_{I}.id_{I}= id.
\end{equation*}
As a result, for each $I$ and each $(k|l) \times(k|l)$ supermatrix $T$, we can see that
\begin{equation*}
T=(T.id_{I}).id_{I}.
\end{equation*}
The following steps may be taken in order to construct the structure sheaf of $_{\nu}Gr$:
\\
Step1: For each $k|l$ multi-index $I$, consider the $\nu$-domain $(V_{I}, \mathcal{O}_{I})$.
\\ 
Step2: Corresponding to each $I$, consider a $(k|l) \times (m|n)$ supermatrix $A^{I}$ which its columns with indices in $I$ together form $id_{I}$. The formal expression $1\nu$ appears when a diagonal entry of $id_{I}$ places in odd part of $A^{I}$.
\\
In addition, the other columns of $A^{I}$, from left to right, and each one from up to down, are filled by even and odd coordinates of $\nu$-domain $\mathcal{O}_{I}$, i.e., 
$x_1, ..., x_k, e_1, ..., e_l, ..., x_{(m- k-1)k+1}, ..., x_{(m- k)k},$ $e_{(m- k-1)l+1}, ..., e_{(m- k)l}, e_{(m- k)l+1} ...e_{(m- k)l+ k}, x_{(m- k)k+1}, ..., x_{(m- k)k+ l}, ...,$
 respectively. Afterwards, each entry, say $a$, that is in a block with opposite parity is replaced by $\nu(a)$.
\\
As an example, consider $_{\nu}Gr_{2|2}(3|3)$ with $I=\{ 1, 2, 3, 6\}$. Then one has
$$\left[ \! \! \! \!
\begin{tabular}{ccccccc}
1   & 0 & 0          & ; & $\nu x_{1}$ & $e_{3}$ & 0\\
0   & 1 & 0          & ; & $\nu x_{2}$ & $e_{4}$ & 0\\
--   &-- & --         & --& --  --            &      --  -- & --\\
0   & 0 & $1\nu$ & ; & $\nu e_{1}$  & $x_{3}$ & 0
\\
0   & 0 & 0          & ; & $\nu e_{2}$ & $x_{4}$ & 1
\end{tabular}
\! \! \! \! \right]$$
The columns of $A^{I}$ with indices in $I$ together form the following supermatrix:
$$M_{I}(A^{I}):= id_{I}=\left[ \! \! \! \!
\begin{tabular}{ccccccc}
1   & 0 & 0          & ; & 0\\
0   & 1 & 0          & ; & 0\\
--   &-- & --         & --& --\\
0   & 0 & $1\nu$ & ; & 0
\\
0   & 0 & 0          & ; & 1
\end{tabular}
\! \! \! \! \right]$$
For each pair multi-indices $I$ and $J$, define the set $V_{IJ}$ to be the largest subset of  $V_{I}$ on which $M_{J}(A^{I}).id_{J}$  is invertible on it.
\\
Step3: On $V_{IJ}$, the equality
\begin{equation*}
\big(M_{J}(A^{I}).id_{J}\big)^{-1}.A^{I} = A^{J}
\end{equation*}
defines a correspondence between even and odd coordinates of $V_{J}$ and rational expressions in $\mathcal{O}_{I}$ appear as corresponding entries of matrices on the two sides of the equality. By (\cite{Varadarajan}, Th 4.3.1), one has a unique homomorphism
\begin{equation*}
\varphi_{IJ}^{*}: \mathcal{O}_{J|_{V_{JI}}} \longrightarrow \mathcal{O}_{I|_{V_{IJ}}}
\end{equation*}
Step4: The homomorphisms $\varphi_{IJ}^{*}$ satisfy the gluing conditions, i.e., for each $I$, $J$ and $K$, we have
\begin{itemize}
\item[1)]
$\varphi_{II}^{*}= id$
\item[2)]
$\varphi_{IJ}^{*} \circ \varphi_{JI}^{*} = id$
\item[3)]
$\varphi_{IK}^{*} \circ \varphi_{KJ}^{*} \circ \varphi_{JI}^{*}= id$
\end{itemize}
In the first condition, $\varphi_{II}^{*}$ is defined by the following equality:
\begin{align*}
\big(M_{I}(A^{I}).id_{I}\big)^{-1}.A^{I} = \big(id_{I}.id_{I}\big)^{-1}.A^{I}
\end{align*}
Since $\big(id_{I}.id_{I}\big)^{-1}= id$, we have $\varphi_{II}^{*}= id$.
\\
For the last condition, note that $\varphi_{KJ}^{*} \circ \varphi_{JI}^{*}$ is obtained from the equality
\begin{equation*}
\Big( M_{I}\Big( \big(M_{J}(A^{K}) .id_{J}\big)^{-1}.A^{K}\Big).id_{I}\Big)^{-1}\Big( \big(M_{J}(A^{K}) .id_{J}\big)^{-1}.A^{K}\Big)= A^{I},
\end{equation*}
For the left hand side of this equality, one has
\begin{align*}
     &\Big(\big(M_{J}(A^{K}) .id_{J}\big)^{-1}. M_{I}(A^{K}).id_{I}\Big)^{-1}\Big( \big(M_{J}(A^{K}) .id_{J}\big)^{-1}.A^{K}\Big)
     \\
   =&\big(M_{I}(A^{K}).id_{I}\big)^{-1} \big(M_{J}(A^{K}) .id_{J}\big)\big(M_{J}(A^{K}) .id_{J}\big)^{-1}. A^{k}
   \\
    =&\big(M_{I}(A^{K}).id_{I}\big)^{-1}.A^{K}= A^{I}
    \end{align*}
Thus the third condition is established.
\\
The second condition may results from other conditions as follows:
\begin{align*}
&\varphi_{IJ}^{*} \circ \varphi_{JI}^{*}=\varphi_{II}^{*}
\\
&\varphi_{II}^{*}=id.
\end{align*}
\section{Super vector bundles}
Here, we recall the definition of super vector bundles and their homomorphisms. Then, we introduce a canonical super vector bundle over $\nu$-grassmannian.
\\
By a \textit{super vector bundle} $\mathcal{E}$ of rank $k|l$ over a supermanifold $(M, \mathcal{O})$, we mean a sheaf of $\mathbb{Z}_{2}$-graded $\mathcal{O}$-modules on $M$ which locally is a free $k|l$ module.
\\
In other words, there exists an open cover $\{U_{\alpha}\}_{\alpha}$ of $M$ such that 
\begin{equation*}
\mathcal{E}(U_{\alpha}) \simeq \big(\mathcal{O}(U_{\alpha}) \big)^{k} \oplus \big(\pi\mathcal{O}(U_{\alpha}) \big)^{l} 
\end{equation*}
or equivalently,
\begin{equation*}
\mathcal{E}(U_{\alpha}) \simeq \mathcal{O}(U_{\alpha}) \otimes_{\mathbb{R}} \big(\mathbb{R}^{k} \oplus \pi(\mathbb{R}^{l}) \big).
\end{equation*}
For example, let $\mathcal{O}_{M}^{k|l} :=(\oplus_{i=1}^{k} \mathcal{O}) \oplus (\oplus_{j=1}^{l} \pi \mathcal{O})$ where $\pi \mathcal{O}$ is an $\mathcal{O}$-module which satisfies
\begin{equation*}
(\pi \mathcal{O})^{ev} = \mathcal{O}^{odd}, \qquad (\pi \mathcal{O})^{odd} = \mathcal{O}^{ev}.
\end{equation*}
The right multiplication is the same as in $\mathcal{O}$ and the left multiplication is as follows:
\begin{equation*}
z(\pi w):= (-1)^{p(z)}\pi(zw) 
\end{equation*}
where $\pi z$ is an element of $\pi\mathcal{O}$.
\\
Hence, $\mathcal{O}_{M}^{k|l}$ is a super vector bundle over the supermanifold $M$.
\\
Let $\mathcal{E}$ and $\mathcal{E}^{\prime}$ be two super vector bundles over a supermanifold $(M, \mathcal{O})$. By a homomorphism from $\mathcal{E}$ to $\mathcal{E}^{\prime}$, we mean an even sheaf homomorphism $\tau: \mathcal{E} \longrightarrow \mathcal{E}^{\prime}$.
\\
Each super vector bundle over $M$ isomorphic to $\mathcal{O}_{M}^{k|l}$, is called a \textit{trivial super vector bundle} of rank $k|l$.
\subsection{Canonical super vector bundle over $\nu$-grassmannian}
Let $I$ be a $k|l$ multi-index and let $(V_{I}, \mathcal{O}_{I})$ be a $\nu$-domain. Consider the trivial super vector bundle
\begin{equation*}
\Gamma_{I}:= \mathcal{O}_{I} \otimes_{\mathbb{R}} \big(\mathbb{R}^{k} \oplus \pi(\mathbb{R}^{l}) \big) = \mathcal{O}_{I} \otimes_{\mathbb{R}} \mathbb{R}^{k|l}
\end{equation*}
By gluing these super vector bundles through suitable homomorphisms, one may construct a super vector bundle $\Gamma$ over $\nu$-grassmannian $_{\nu}Gr$. For this, consider a basis $\{e_{1}, \cdots, e_{k}, f_{1}, \cdots, f_{l}\}$ for $\mathbb{R}^{k|l}$ and set
\begin{equation*}
m:= \big(M_{J}(A^{I}) .id_{J}\big)^{-1}
\end{equation*}
where $\big(M_{J}(A^{I}) .id_{J}\big)^{-1}$ is introduced in subsection 2.2. Gluing morphisms are defined as follows:
\begin{equation*}
\psi_{IJ}^{*}: \Gamma_{J|_{V_{JI}}} \longrightarrow \Gamma_{I|_{V_{IJ}}}
\end{equation*}
\begin{equation*}
a \otimes e_{i} \longmapsto \varphi_{IJ|_{V_{JI}}}^{*}(a) \Big( \sum_{t \leq k} m_{it} \otimes e_{t} + \sum_{t > k}m_{it}\otimes f_{t}\Big),
\end{equation*}
where the elements $m_{it}$ are the entries of the $i$-th column of the supermatrix $m$. The morphisms $\psi_{IJ}^{*}$ satisfy the gluing conditions. So $\Gamma_{I}^,$s may glued together to form a super vector bundle denoted by $\Gamma$.
\section{Gauss supermaps}
In common geometry, a Gauss map is defined as a map from the total space of a vector bundle, say $\xi$, to a Euclidean space such that its restriction to any fiber is a monomorphism. Equivalently, one may consider a $1-1$ strong bundle map from $\xi$ to a trivial vecor bundle. The Gauss map induces a homomorphism between the vector bundle and the canonical vector bundle on a grassmannian $Gr_{k}(n)$ with sufficiently large value of $n$. A simple method for making such a map is the use of coordinate representation for $\xi$. In this section, for constructing a Gauss supermap of a super vector bundle, one may use the same method.
\\
We call a super vector bundle $\mathcal{E}$ over a supermanifold $(M, \mathcal{O})$ is of finite type, whenever there is a finite open cover $\{U_{\alpha}\}_{\alpha=1}^{t}$ for $M$ such that the restriction of $\mathcal{E}$ to $U_{\alpha}$ is trivial, i.e., there exists isomorphisms
\begin{equation*}
\psi_{\alpha}^*: \mathcal{E}|_{U_{\alpha}} \overset{\simeq}{\longrightarrow} \mathcal{O}|_{U_{\alpha}} \otimes _{\mathbb{R}} \mathbb{R}^{k|l}.
\end{equation*}
A Gauss supermap over $\mathcal{E}$ is a homomorphism from $\mathcal{E}$ to the trivial super vector bundle over $(M, \mathcal{O})$ so that its kernel is trivial.
\\
Let $\{e_{1}, \cdots, e_{k}, f_{1}, \cdots, f_{l}\}$ be a basis for $\mathbb{R}^{k|l}$ so that $\{e_{i}\}$ and $\{f_{j}\}$ are respectively bases for $\mathbb{R}^{k}$ and $\pi(\mathbb{R}^{l})$, then $B:=\{1 \otimes e_{1}, \cdots, 1 \otimes e_{k}, 1 \otimes f_{1}, \cdots, 1 \otimes f_{l}\}$ is a generator for the $\mathcal{O}(U_{\alpha})$-module, $\mathcal{O}(U_{\alpha}) \otimes_{\mathbb{R}} \mathbb{R}^{k|l}$.
\\
Set
\begin{equation}\label{salpha}
s_{i}^{\alpha}:= \psi_{\alpha}^{*^{-1}} (1 \otimes e_{i}), \qquad t_{j}^{\alpha}:= \psi_{\alpha}^{*^{-1}} (1 \otimes f_{j}).
\end{equation}
So $(\psi_{\alpha}^{*})^{-1} (B)$ is a generator for $\mathcal{E}(U_{\alpha})$ as an $\mathcal{O}(U_{\alpha})$-module.
\\
Choose a partition of unity $\{\rho_{\alpha}\}_{\alpha=1}^{t}$ subordinate to the covering $\{U_{\alpha}\}_{\alpha=1}^{t}$. Considering $s$ as a global section of $\mathcal{E}(M)$, we can write 
\begin{equation*}
s=\sum_{\alpha=1}^{t} \rho_{\alpha}r_{\alpha}(s)
\end{equation*}
where $r_{\alpha}$ is the restriction morphism. In addition, one has
\begin{equation*}
r_{\alpha}(s)=\sum_{i=1}^{k} \lambda_{i}^{\alpha}s_{i}^{\alpha} + \sum_{j=1}^{l} \delta_{j}^{\alpha}t_{j}^{\alpha}, \qquad \lambda_{i}^{\alpha}, \delta_{j}^{\alpha} \in \mathcal{O}(U_{\alpha}).
\end{equation*}
By the last two equalities, we have
\begin{equation*}
s=\sum_{\alpha=1}^{t} \rho_{\alpha} \Big( \sum_{i=1}^{k} \lambda_{i}^{\alpha}s_{i}^{\alpha} + \sum_{j=1}^{l} \delta_{j}^{\alpha}t_{j}^{\alpha}\Big) = \sum_{\alpha=1}^{t}\sum_{i=1}^{k} \sqrt{\rho_{\alpha}}\lambda_{i}^{\alpha}.
\sqrt{\rho_{\alpha}}s_{i}^{\alpha}+ \sum_{\alpha=1}^{t}\sum_{j=1}^{l} \sqrt{\rho_{\alpha}}\delta_{j}^{\alpha}.\sqrt{\rho_{\alpha}}t_{j}^{\alpha},
\end{equation*}
where $\sqrt{\rho_{\alpha}}s_{i}^{\alpha}$ and $\sqrt{\rho_{\alpha}}t_{j}^{\alpha}$ are even and odd sections of $\mathcal{E}(M)$ repectively, and $\sqrt{\rho_{\alpha}}\lambda_{i}^{\alpha}$ and $\sqrt{\rho_{\alpha}}\delta_{j}^{\alpha}$
are sections of $\mathcal{O}(M)$. So $A:= \{\sqrt{\rho_{\alpha}}s_{i}^{\alpha}\}_{\alpha, i} \cup \{\sqrt{\rho_{\alpha}}t_{j}^{\alpha}\}_{\alpha, j}$ is a generating set of $\mathcal{E}(M)$.
\\
Now, for each $\alpha$, consider the following monomorphism between $\mathcal{O}(U_{\alpha})$-modules:
\begin{equation*}
i_{\alpha}: \mathcal{O}(U_{\alpha}) \otimes_{\mathbb{R}} \mathbb{R}^{k|l} \longrightarrow \mathcal{O}(U_{\alpha}) \otimes_{\mathbb{R}} \mathbb{R}^{tk|tl}
\end{equation*}
\begin{equation*}
1 \otimes e_{i} \longmapsto 1 \otimes e_{(\alpha-1)k+i}
\end{equation*}
\begin{equation*}
1 \otimes f_{j} \longmapsto 1 \otimes f_{(\alpha-1)l+j}.
\end{equation*}
Set
\begin{equation}
g(s):= \sum_{\alpha=1}^{t} \rho_{\alpha}. i_{\alpha} \circ \psi_{\alpha}^* \circ r_{\alpha}(s).
\end{equation}
It is easy to see that $g$ is a Gauss supermap of $\mathcal{E}(M)$.
\subsection{Gauss supermatrix}
Now, we are going to obtain the matrix of the gauss supermap $g$.
\\
By a Gauss supermatix associated to super vector bunle $\mathcal{E}$, we mean a supermatrix, say $G$, which is obtained as follows with respect to the generating set $A$:
\begin{equation*}
g\big(\sqrt{\rho_{\beta}}s_{j}^{\beta}\big) = \sum_{\alpha=1}^{t} \rho_{\alpha}. i_{\alpha} \circ \psi_{\alpha}^* \circ r_{\alpha}\big( \sqrt{\rho_{\beta}}s_{j}^{\beta}\big)
\end{equation*}
where $g$ is a Gauss supermap over $\mathcal{E}$.
\\
By \eqref{salpha}, we have
\begin{align}\label{grho}
g\big(\sqrt{\rho_{\beta}}s_{j}^{\beta}\big) &= \sum_{\alpha=1}^{t} \rho_{\alpha}. i_{\alpha} \circ \psi_{\alpha}^* \circ \psi_{\beta}^{*^{-1}}\big( \sqrt{\rho_{\beta}}e_{j}\big) \nonumber
\\
                                                                 &=\sum_{i=1}^{k}\sum_{\alpha=1}^{t} \rho_{\alpha} \sqrt{\rho_{\beta}} a_{ij}^{\alpha\beta}e_{(\alpha-1)k+i}+\sum_{i=1}^{l}\sum_{\alpha=1}^{t}\rho_{\alpha} \sqrt{\rho_{\beta}}a_{(k+i)j}^{\alpha\beta}f_{(\alpha-1)l+i} 
\end{align}
where $[a_{ij}^{\alpha\beta}]$ is a matrix of $\psi_{\alpha}^* \circ \psi_{\beta}^{*^{-1}}$ relative to the basis $B$. The natural ordering on $\{i_{\alpha}e_{i}\}_{\alpha, i}$ and $\{i_{\alpha}f_{s}\}_{\alpha, s}$ induces an ordering on their coefficients in \eqref{grho}. Let G be a $tk|tl \times tk|tl$ standard supermatrix. Fill the even and odd top blocks of $G$ by these coefficients according to their parity from left to right along the $\big((\beta-1)k+j\big)$-th row, $1 \leq j \leq k$, $1 \leq \beta \leq t$. Similarly, by coefficients in the decomposition of $g\big(\sqrt{\rho_{\beta}}t_{r}^{\beta}\big)$, one may fill the odd and even down blocks of $G$ along the $\big((\beta -1)k+r\big)$-th row, $1 \leq r \leq l$, $1 \leq \beta \leq t$.
\begin{example}
For $k|l=2|1$ and a suitable covering with two elements, i.e., $t=2$, we have
\begin{equation*}
g\big(\sqrt{\rho_{2}}s_{2}^{2}\big)= \rho_{1} \sqrt{\rho_{2}} a_{12}^{12}e_{1} +\rho_{1} \sqrt{\rho_{2}} a_{22}^{12}e_{2}+\rho_{2} \sqrt{\rho_{2}} a_{12}^{22}e_{3}+\rho_{2} \sqrt{\rho_{2}} a_{22}^{22}e_{4}+\rho_{1} \sqrt{\rho_{2}} a_{32}^{12}f_{1}+\rho_{2} \sqrt{\rho_{2}} a_{32}^{22}f_{2}  
\end{equation*}
Then the $4$-th row of the associated supermatrix $G$ is as below:
$$\left[ \! \! \! \!
\begin{tabular}{ccccccc}
                                                              &                                                          &                                                           &                                                                                               \vdots 
                                                              & ;                          &                                                            &
                                             \\
$\rho_{1} \sqrt{\rho_{2}} a_{12}^{12}$   & $\rho_{1} \sqrt{\rho_{2}} a_{22}^{12}$ & $\rho_{2} \sqrt{\rho_{2}} a_{12}^{22}$ & $\rho_{2} \sqrt{\rho_{2}} a_{22}^{22}$  & ;                         &$\rho_{1} \sqrt{\rho_{2}} a_{32}^{12}$ & $\rho_{2} \sqrt{\rho_{2}} a_{32}^{22}$\\
-- --         -- --            -- --                        &     -- --       -- --    -- --                         & -- --     -- --     -- --                              & -- --   -- --  -- --                                                         & -- --   -- --   -- --    &-- --  -- -- -- --                                        &   -- --   -- --   -- --  \\
                                                               &                                                          &                                                           &                                                                                               \vdots 
                                                               & ;  &                                                             &
                                             \\
\end{tabular}
\! \! \! \! \right].$$
\end{example}
On the other hand, one may consider a covering $\{U_{\alpha}\}_{\alpha}$ so that for each $\alpha$, we have an isomorphism
\begin{equation}\label{isomorphism}
\mathcal{O}(U_{\alpha}) \overset{\simeq}{\longrightarrow} \mathbf{C}^{\infty}(\mathbb{R}^{m}) \otimes _{\mathbb{R}} \wedge \mathbb{R}^{n}
\end{equation}
Let $\nu$ be an odd involution on $\mathbf{C}^{\infty}\big(\mathbb{R}^{(k^{2}+l^{2})(t-1)}\big) \otimes _{\mathbb{R}} \wedge \mathbb{R}^{2kl(t-1)}$ preserving $\mathbf{C}^{\infty}(\mathbb{R}^{m}) \otimes _{\mathbb{R}} \wedge \mathbb{R}^{n}$ as a subalgebra. Thus, it induces an odd involution on $\mathcal{O}(U_{\alpha})$ through the isomorphism \eqref{isomorphism} which is denoted by the same notation $\nu$.
\begin{theorem}
Let $\cal{E}$ be a super vector bundle over a supermanifold $(M, \cal{O})$ and let $G$ be a Gauss supermatrix associated to $\cal{E}$. Then the Gauss supermatrix induces a homomorphism from $\mathcal{G}$, the structure sheaf of   $_{\nu}Gr$, to $\mathcal{O}$.
\end{theorem}
\textit{Proof}.
Let $h$ be an element of $\mathcal{G}(M)$, and $\{\rho_{I}^{\prime}\}$ be a partition of unity subordinate to the covering $\{V_{I}\}_{I \subseteq \{1, \cdots, tk\}}$, then one has
\begin{equation}
h= \sum_{\substack{I \subseteq \{1, \cdots, tk\},\\ |I|=k+l}} \rho_{I}^{\prime}.h|_{V_{I}}
\end{equation}
Consider the rows of $G$ with indices in $I$ as a $(k+l) \times t(k+l)$ supermatrix and name it $G(I)$. Then multiply it by $id_{I}$ from left, i.e., $id_{I}. G(I)$ and delete the columns with indices in $I$, we get
$$\left[ \! \! \! \!
\begin{tabular}{ccccccc}
$y_{11}^{I}$  $\qquad$     & $\cdots$    $\qquad$       &$y_{1\big((t-1)(k+l)\big)}^{I}$\\
$\vdots$          $\qquad$    &  $\ddots$   $\qquad$        &$\vdots$ \\
$y_{(k+l)1}^{I}$ $\qquad$ &  $\cdots$   $\qquad$       &$y_{(k+l)\big((t-1)(k+l)\big)}^{I}$
\end{tabular}
\! \! \! \! \right].$$
Note that all entries of this supermatrix are sections of $\mathcal{O}(M)$.
\\
Let $A^{I}$  be the matrix introduced in subsection 2.2 and let $x_{ij}^{I}$ be its entry out of $M_{I}(A^{I})=id_{I}$. Then the correspondence $x_{ij}^{I} \longmapsto y_{ij}^{I}$ defines a homomorphism 
\begin{equation*}
\varphi_{I}^{*}: \mathcal{O}_{I}(V_{I}) \longrightarrow \mathcal{O}(M).
\end{equation*}
Now, for each global section $h$ of $\mathcal{G}$, one may define 
\begin{equation*}
\widetilde{h} = \sum_{\substack{I \subseteq \{1, \cdots, tk\},\\ |I|=k+l}} \varphi_{I}^{*}( \rho_{I}^{\prime}.h|_{V_{I}})
\end{equation*}
Then the correspondence 
\begin{equation}\label{sigma}
\sigma^*: h \longmapsto \widetilde{h}
\end{equation}
is a  well-defined homomorphism from $\mathcal{G}(Gr_{k}^{tk} \times Gr_{l}^{tl})$ to $\mathcal{O}(M)$ and so induces a smooth map $\widetilde{\sigma}$, from $M$ to $Gr_{k}^{tk} \times Gr_{l}^{tl}$ \cite{Roshandel}.
\begin{flushright}
$\square$
\end{flushright}
The homomorphism $\sigma=(\widetilde{\sigma},\sigma^*)$ is called the associated morphism with the Gauss supermap $g$.
\section{Pullback of the canonical super vector bundle}
\begin{theorem}
The super vector bundle $\mathcal{E}$ and the pullback of $\Gamma$ ,the canonical super vector bundle over $_{\nu}Gr$, under $\sigma$ are isomorphic.
\end{theorem}
\textit{Proof}.
First note that one may define a $\mathcal{G}$-module structure on $\mathcal{O}(M)$ by applying $\sigma$ in \eqref{sigma} as follows:
\begin{equation*}
a * b:= \sigma(a).b, \qquad a \in \mathcal{G}(Gr_{k}^{tk} \times Gr_{l}^{tl}), \qquad b \in \mathcal{O}(M).
\end{equation*}
The pullback of $\Gamma$ along $\sigma$ is defined as $\mathcal{O} \otimes _{\mathcal{G}}^{\sigma} \Gamma$.
\\
We show that this module is isomorphic to $\mathcal{E}$. Let $s^{\prime}$ be a global section on $\Gamma$. One has
\begin{equation*}
s^{\prime}= \sum_{\substack{I \subseteq \{1, \cdots, tk\},\\ |I|=k+l}} \rho_{I}^{\prime}. r_{I}^{\prime}(s^{\prime})
\end{equation*}
where $\{\rho_{I}^{\prime}\}$ is the partition of unity of $_{\nu}Gr$ subordinate to the open cover $\{V_{I}\}$, and $r_{I}^{\prime}$ is the restriction morphism giving sections over $V_{I}$. On the other hand, one may write each section $r_{I}^{\prime}(s^{\prime})$ as below:
\begin{equation*}
r_{I}^{\prime}(s^{\prime}) = \sum_{j=1}^{k+l} h_{j}^{I} s_{j}^{\prime^{I}}
\end{equation*}
where $s_{j}^{\prime^{I}}$ are generators of $\Gamma(V_{I})$ and the coefficients $h_{j}^{I}$ are the sections of $\mathcal{O}_{I}$. Therefore, we can write 
\begin{equation*}
s^{\prime}= \sum_{\substack{I \subseteq \{1, \cdots, tk\},\\ |I|=k+l}} \sum_{j=1}^{k+l} (\rho_{I}^{\prime}h_{j}^{I})s_{j}^{\prime^{I}}
\end{equation*}
Note that each row of $G$ is in correspondence with a section in the generator set $A$. So there is a morphism from the pullback of $\Gamma$ to $\mathcal{E}$ as
\begin{align}\label{sigmatensor}
&\mathcal{O} \otimes_{\mathcal{G}}^{\sigma} \Gamma \longrightarrow \mathcal{E} \nonumber
\\
&u \otimes s^{\prime} \longmapsto u.\delta(s^{\prime})
\end{align}
where $\delta(s^{\prime})$ is 
\begin{equation*}
\sum_{\substack{I \subseteq \{1, \cdots, tk\},\\ |I|=k+l}} \sum_{\substack{j=1}}^{k+l} \sigma(\rho_{I}^{\prime}h_{j}^{I}). s_{j}^{I}
\end{equation*}
and $s_{j}^{I}$ is the section corresponding to the $j$-th row of $G(I)$ (cf. subsection 4.1).
\\
One may show that the morphism in \eqref{sigmatensor} is an isomorphism. To this end, first note that every locally isomorphism between two sheaves of $\mathcal{O}$-modules with the same rank is a globally isomorphism. Also for the super vector bundle $\Gamma$ of rank $k|l$ over $\mathcal{G}$, one can write a locally isomorphism
\begin{equation*}
\mathcal{O} \otimes_{\mathcal{G}}^{\sigma} \Gamma \overset{\simeq}{\longrightarrow} \mathcal{O} \otimes_{\mathbb{R}}\mathbb{R}^{k|l},
\end{equation*}
because for each sufficiently small open set $V$ in $Gr_{k}^{tk} \times Gr_{l}^{tl}$ one can write
\begin{equation*}
\Gamma(V) \simeq \mathcal{G}(V) \otimes_{\mathbb{R}}\mathbb{R}^{k|l} 
\end{equation*}
and then
\begin{equation*}
\mathcal{O} \big(\widetilde{\sigma}^{-1}(V)\big) \otimes_{\mathcal{G}}^{\sigma} \Gamma(V) \simeq \mathcal{O}\big(\widetilde{\sigma}^{-1}(V)\big) \otimes_{\mathcal{G}}^{\sigma} \mathcal{G}(V) \otimes_{\mathbb{R}} \mathbb{R}^{k|l}
\end{equation*}
This shows that the morphism in \eqref{sigmatensor} may be represented locally by the following isomorphism:
\begin{equation}
\mathcal{O} \otimes \mathcal{G}\otimes \mathbb{R}^{k|l} \to \mathcal{O}\otimes \mathbb{R}^{k|l}.
\end{equation}
Thus \eqref{sigmatensor} defines a global isomorphism. 
\begin{flushright}
$\square$
\end{flushright}
\section{Homotopy properties of Gauss supermaps and their associated morphisms}

Let $\mathcal{O}_{M}^{m|n}$ and $\mathcal{O}_{M}^{m^{\prime}|n^{\prime}}$ be two trivial super vector bundles where $m^{\prime}=2m-k$, $n^{\prime}=2n-l$; then, one can write the inclusion homomorphisms
\begin{equation*}
J^{e}, J^{o}, J: \mathcal{O}(M) \otimes_{\mathbb{R}} \mathbb{R}^{m|n} \longrightarrow \mathcal{O}(M) \otimes_{\mathbb{R}} \mathbb{R}^{2m|2n}
\end{equation*}
by the conditions
\begin{align*}
J^{e}: &1 \otimes e_{i} \longmapsto 1 \otimes e_{2i}               &J^{o}: 1 \otimes e_{i}  \longmapsto 1 \otimes e_{2i-1} &\qquad \quad J: 1 \otimes e_{i} \longmapsto 1 \otimes e_{i} \\
           &1 \otimes f_{j} \longmapsto 1 \otimes f_{2j}                & 1 \otimes f_{j} \longmapsto 1 \otimes f_{2j-1}             &\qquad \qquad \quad 1 \otimes f_{j} \longmapsto 1 \otimes f_{j}
\end{align*}
Now, let $(V_I, \mathcal{O}_I)$ be $\nu$-domains introduced in subsection 2.2. In addition assume $(W_J, \mathcal{O}_J)$ be $\nu$-domains of dimension $2p|2q$. For each $k|l$ multi-index $I=\{i_{1}, \cdots, i_{k+l}\} \subset \{1, ..., m+n\}$, one can associate the following multi-indices
\begin{equation*}
I^{e}:=\{2i_{1}, \cdots, 2i_{k+l}\},
\end{equation*}
\begin{equation*}
I^{o}:=\{2i_{1}-1, \cdots, 2i_{k+l}-1\},
\end{equation*}
\begin{equation*}
\bar{I}:=\{i_{1}, \cdots, i_{b}, i_{b+1}+m-k, \cdots, i_{k+l}+m-k\},
\end{equation*}
where $i_a \in I$ ,$1 \leq i_a \leq k+l$, and $i_b$ is an element of $I$ for which $i_b \leq m \leq i_{b+ 1}$.
\\
So the maps $J^{e}$, $J^{o}$ and $J$ induce the homomorphisms
\begin{equation*}
\bar{J}^{e}, \bar{J}^{o}, \bar{J}: _{\nu}Gr(m|n) \longrightarrow _{\nu}Gr(m^{\prime}|n^{\prime}).
\end{equation*}
In fact, $(\bar{J}^{e})^*|_{W_{I^{e}}}$ is obtained by
\begin{align*}
&\mathcal{O}_{I^{e}}(W_{I^{e}}) \longrightarrow \mathcal{O}_{I}(V_{I})
\\
&\left\{
\begin{array}{rl}
y_{i(2j-1)} \longmapsto x_{ij}, \quad i=1, \cdots, k+l, \quad j=1, \cdots, m+n-k-l
\\
\text{other generators} \longmapsto 0 \hspace{6cm}
\end{array}\right.
\end{align*}
\begin{theorem}
Let $f, f_{1}:(M, \mathcal{O}) \longrightarrow _{\nu}Gr(m|n)$ are induced by the gauss supermaps $g$ and $g_{1}$. Then, $\bar{J}f$ and $\bar{J}f_{1}$ induced by $Jg$ and $Jg_{1}$ are homotopic.
\end{theorem}
\textit{Proof}.
Consider the homomorphisms $J^{e}g$ and $J^{o}g_{1}$, with the induced homomorphisms $\bar{J}^{e}f$ and $\bar{J}^{o}f$. One can define a family of homomorphisms
\begin{align*}
&F_{t}: \mathcal{E}(M) \longrightarrow \mathcal{O} \otimes_{\mathbb{R}} \mathbb{R}^{2m|2n}
\\
&\varphi \longmapsto (1-t).(J^{e}g)(\varphi)+t.(J^{o}g_{1})(\varphi)
\end{align*}
where $F_{0}=J^{e}g$ and $F_{1}=J^{o}g_{1}$. By section 4, a family of morphisms $\bar{F}_{t}$ from $(M,\mathcal{O})$ to $_{\nu}Gr(m^{\prime}|n^{\prime})$ are induced. Obviousely  $\bar{F}_{0}=\bar{J}^{e}f$ and $\bar{F}_{1}=\bar{J}^{o}f_{1}$, thus $\bar{F}_{t}$ is a homotopy from $\bar{J}f$ to $\bar{J}f_1$.

\end{document}